\documentclass {amsart}
\usepackage[pdftex]{graphicx}
\usepackage[pdftex]{color}
\usepackage{graphbox}
\usepackage{bigints}
\usepackage{scalerel}
\def\stretchint#1{\vcenter{\hbox{\stretchto[440]{\displaystyle\int}{#1}}}}

\newcommand{\Df}{{\hc \;\partial\!} f}
\newcommand{\ds}{\displaystyle}

\newcommand{\ve}{\varepsilon}
\newcommand{\hc}{}

\makeatletter
\def\Ddots{\mathinner{\mkern1mu\raise\p@
\vbox{\kern7\p@\hbox{.}}\mkern2mu
\raise4\p@\hbox{.}\mkern2mu\raise7\p@\hbox{.}\mkern1mu}}
\makeatother

\title{Flat Map of a Sphere via Stress Minimization}
\author{Robert J. Vanderbei}
	\thanks{Department of Operations Research and Financial Engineering,
	Princeton University,
	Princeton, NJ 08544 ({\tt rvdb@princeton.edu})
	}

\begin{document}

\begin{figure}
\includegraphics[width=5in]{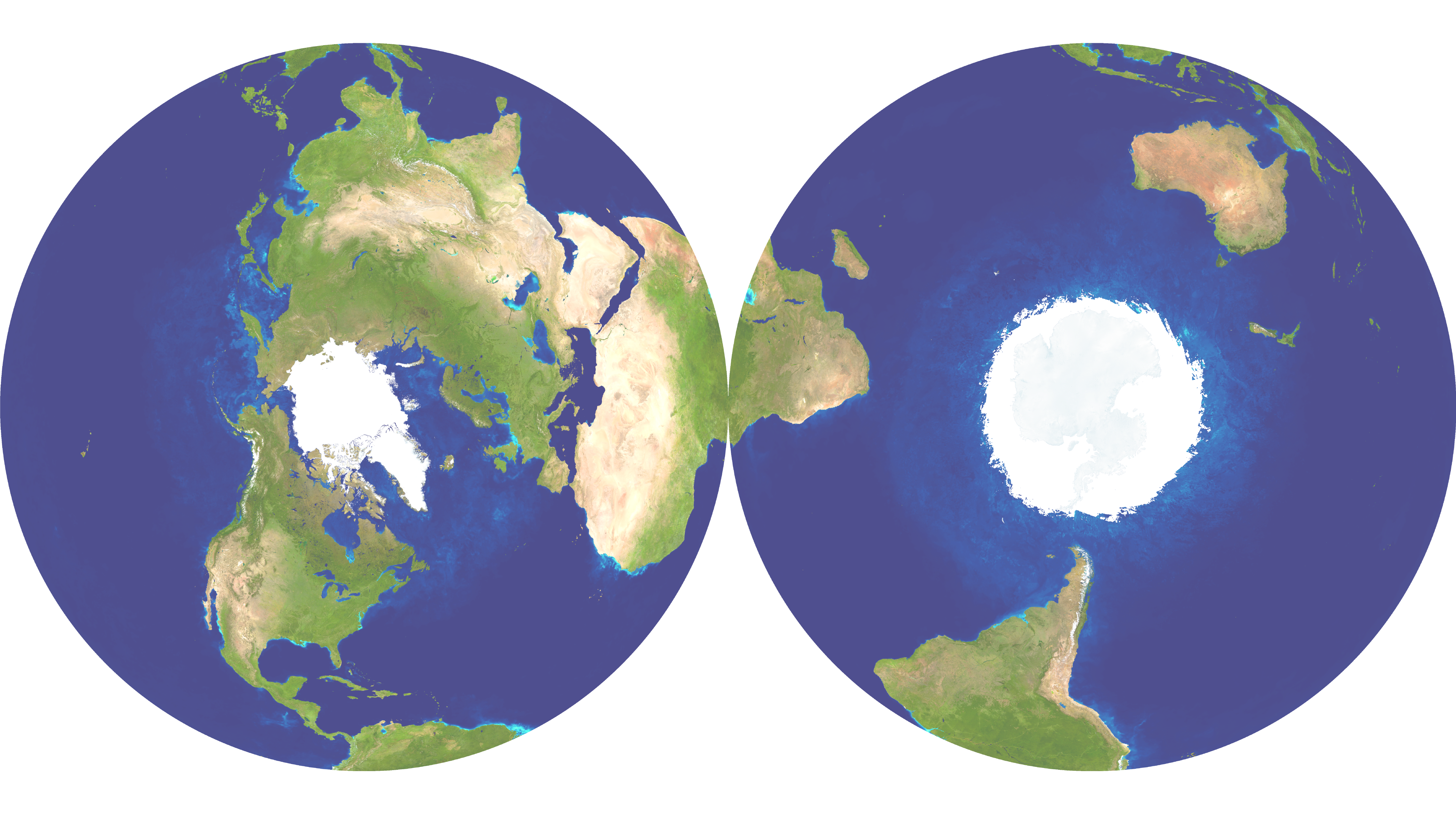}
\caption{The two sides of the flat disk map of the Earth.}
\label{fig3}
\end{figure}

\begin{abstract}
In this paper we describe a mathematically interesting but relatively minor
improvement to the Gott-Goldberg-Vanderbei (GGV) map projection.
This new projection can be described as what one would get by making a spherical
rubber ball representation of the Earth and then stretching the ball circularly
around the equator until the Northern and Southern hemispheres flatten to a
disk.  It is interesting that this new projection is very similar to but not
exactly the same as the GGV projection.  And, the mathematics required to solve
this flattening problem is a very nice example of using the calculus of
variations to solve an infinite dimensional optimization problem.
\end{abstract}

\maketitle


\section{Introduction}

For some thousands of years it has been understood that
our planet Earth is not a flat object -- it's shape is spherical.  
The easiest way to prove this is to observe
multiple partial lunar eclipses.  Lunar eclipses take place when the moon is fully
illuminated, i.e., during a full moon.  And, at time of totality, the Sun and
Moon are in opposite directions from our perspective here on Earth.  A full
Moon rises in the east just after sunset, i.e., when the
Sun is just below the horizon in the west.   Six hours later, the full Moon is
directly overhead and the Sun is straight below us.   So,
when the Moon is full, the Earth is roughly between the Sun and the Moon.   And then,
sometimes during this full phase, we see a shadow being cast on the Moon.
It's not hard to understand that that's Earth's shadow.   And, even though we
don't see the entire shadow of the Earth (because the Earth is bigger than the
Moon), it's pretty obvious that what we see is circular in shape.
Now, some have argued that if the Earth is a flat disk it also would make a
round circular shadow.  That would be true at midnight, when the partially
eclipsed full Moon is roughly straight overhead.   But, if the eclipse were to
take place shortly after sunset (or shortly before sunrise), then the shadow of
the disk would be at a sharp angle to the perpendicular direction and therefore
the shadow would be very elliptical in shape (or maybe cylindrical if the depth
below the disk is significant).  The shape of the eclipse never
appears elliptical.  It's always nice and circular in appearance and hence the
Earth is a sphere.

\begin{figure}
\begin{minipage}{5in}
\centering
\includegraphics[width=1.6in,align=c]{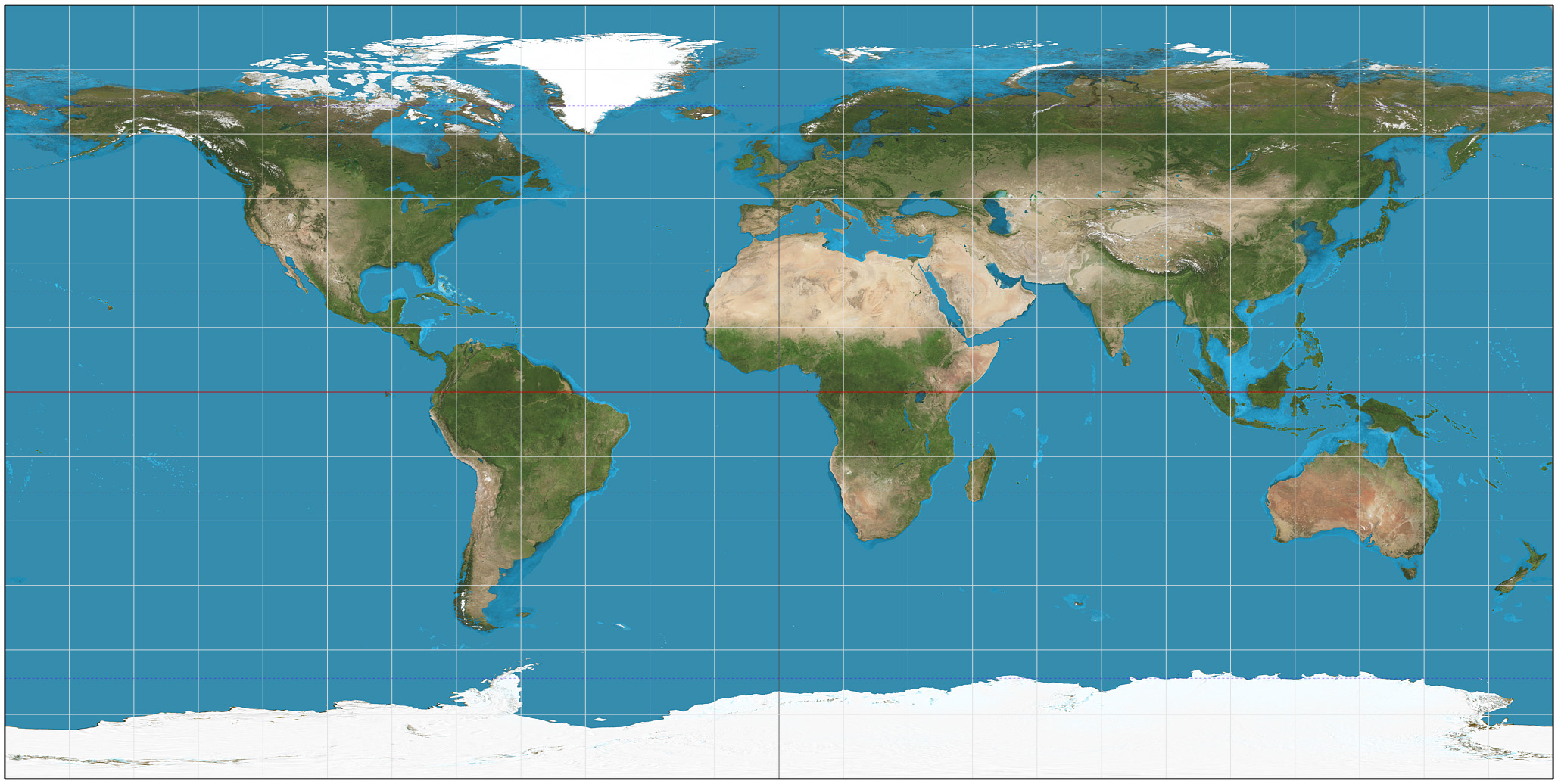}
\includegraphics[width=1.6in,align=c]{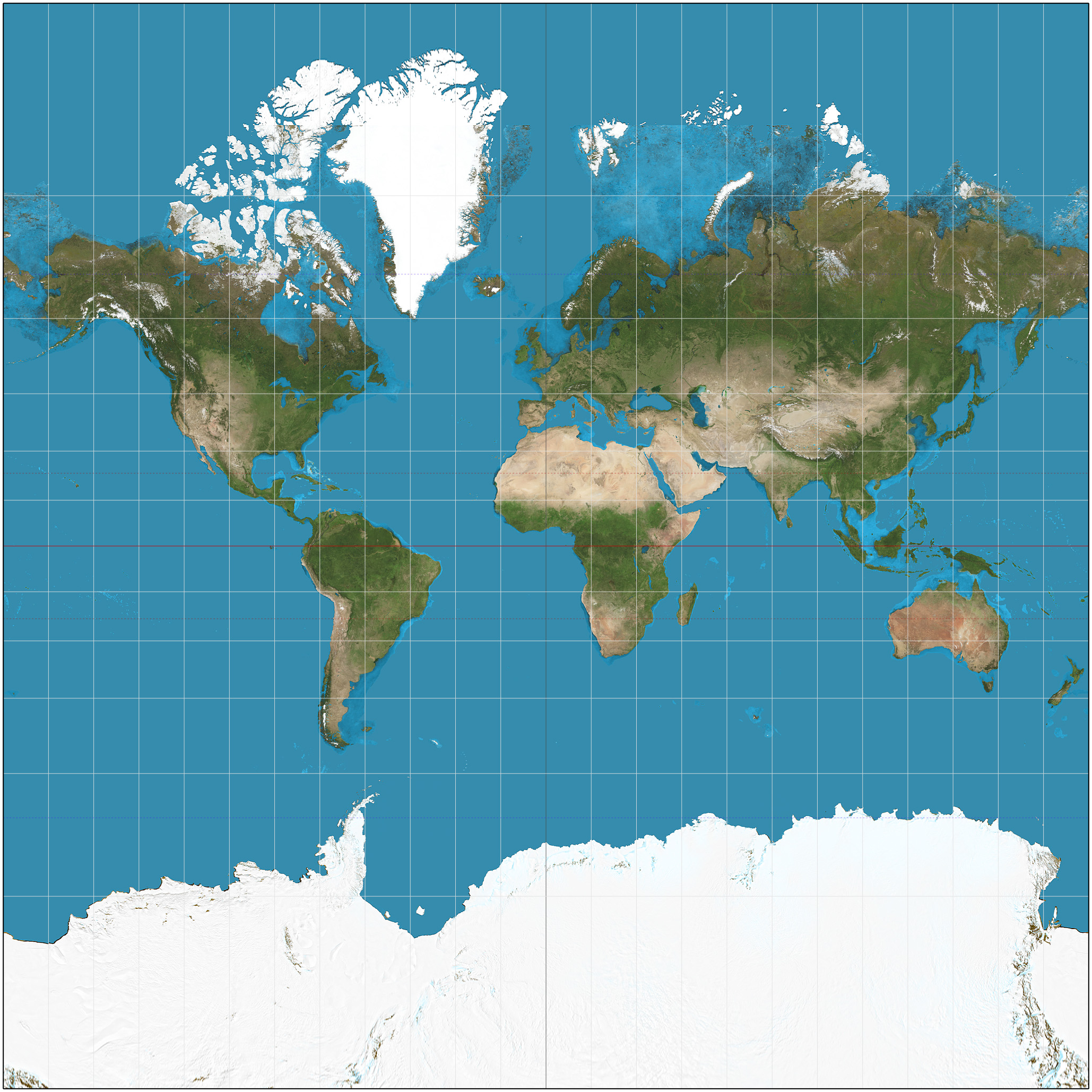} 
\includegraphics[width=1.6in,align=c]{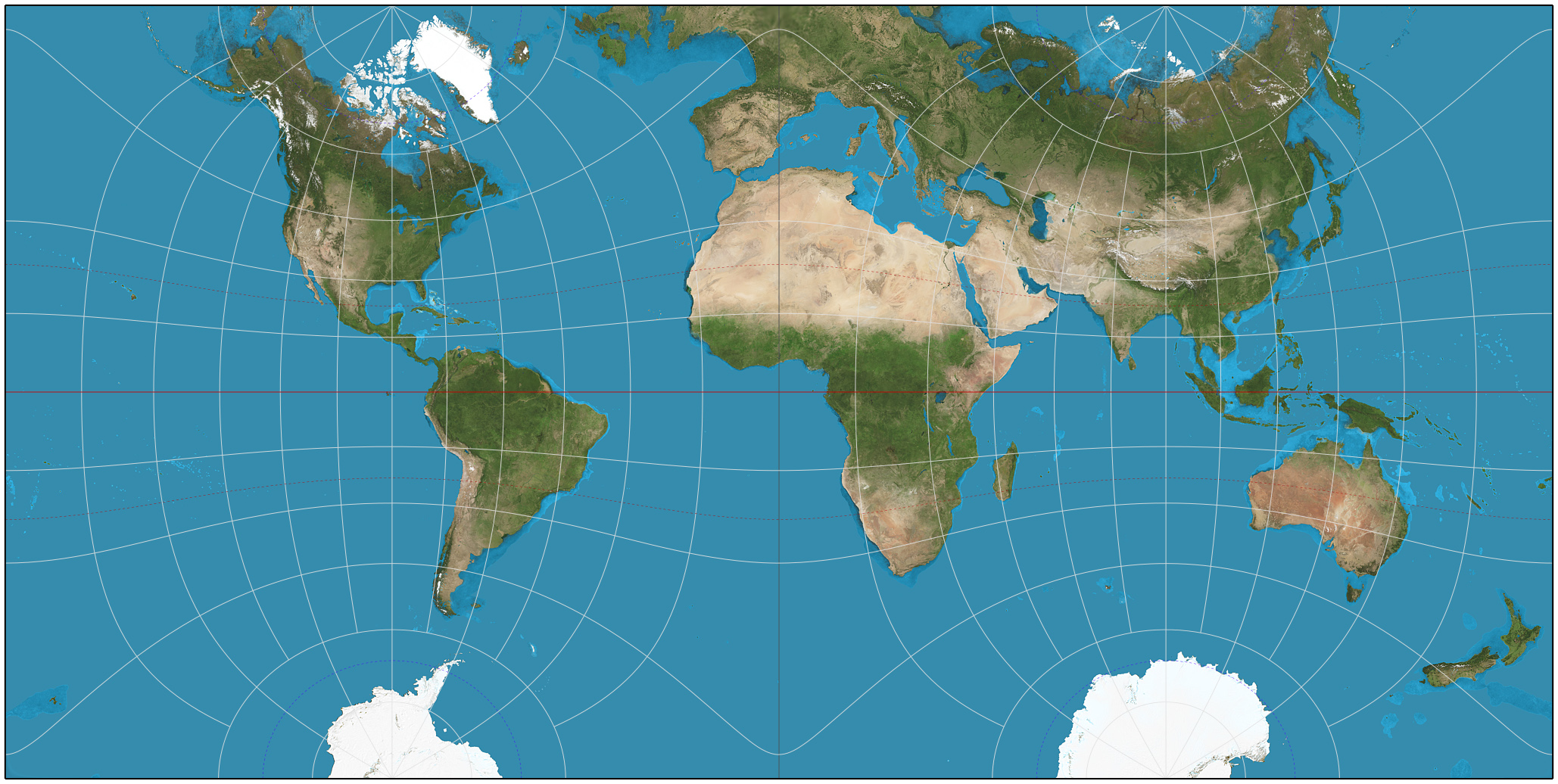}
\\
\includegraphics[width=1.6in,align=c]{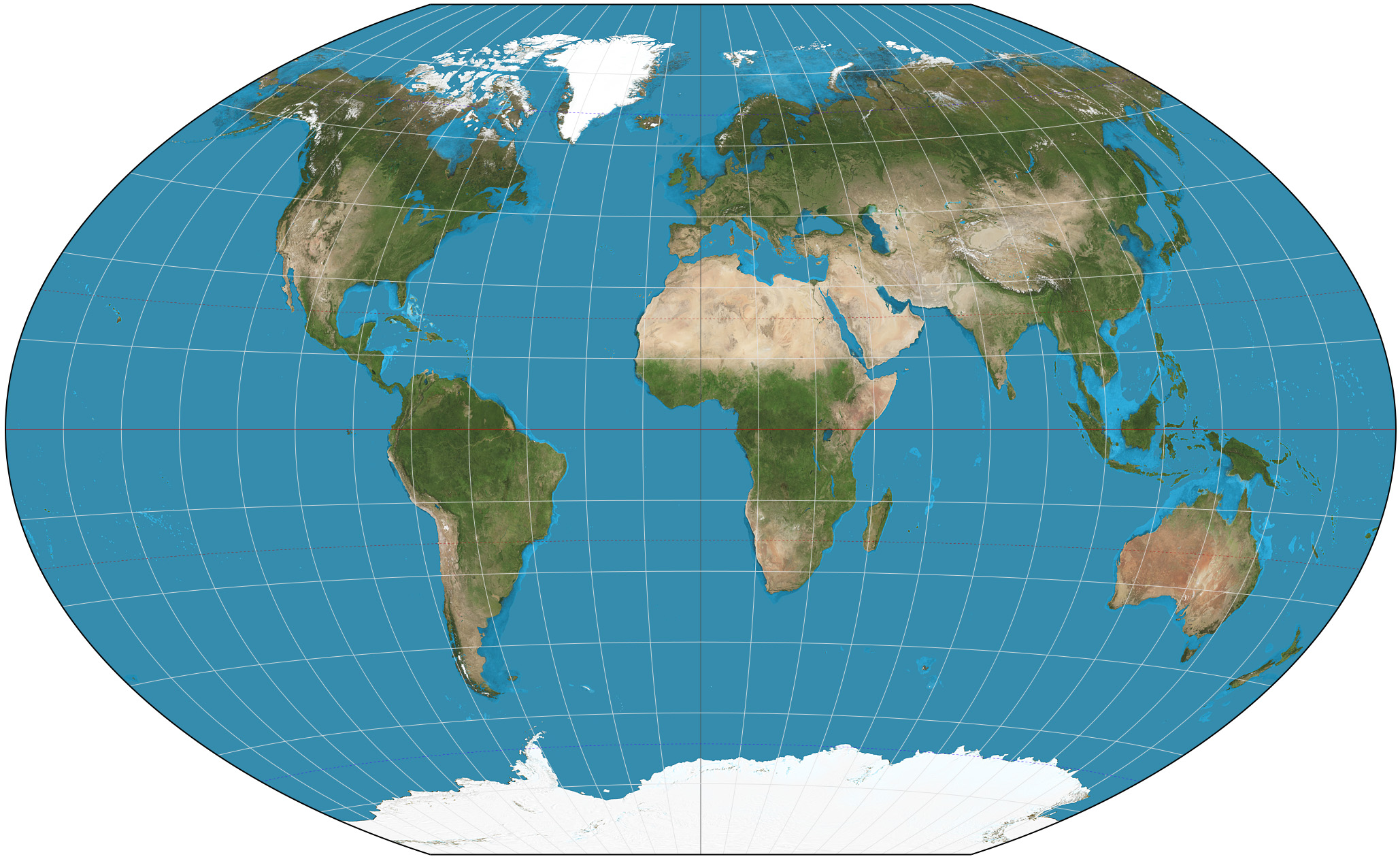}
\includegraphics[width=1.6in,align=c]{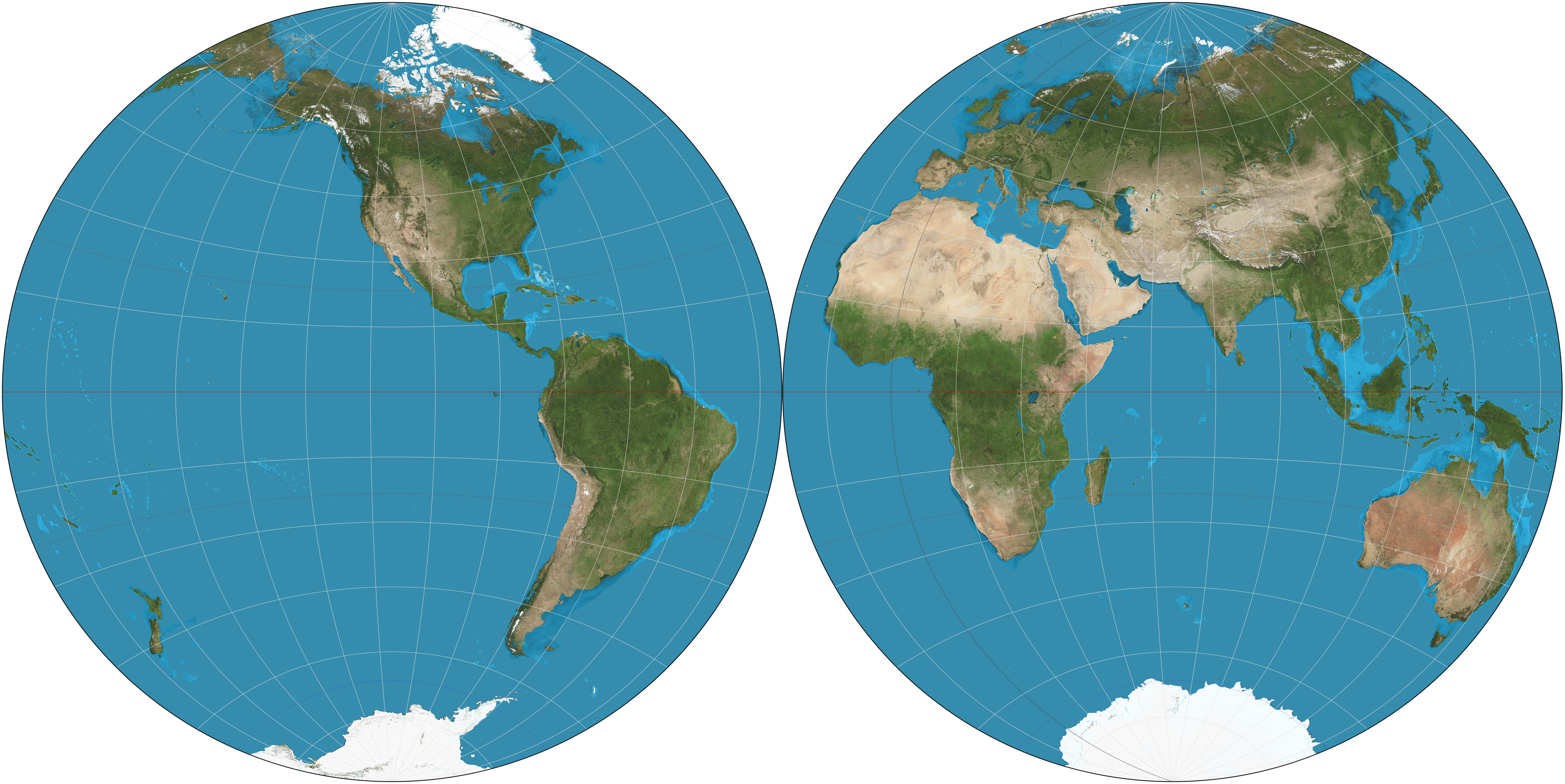}
\includegraphics[width=1.6in,align=c]{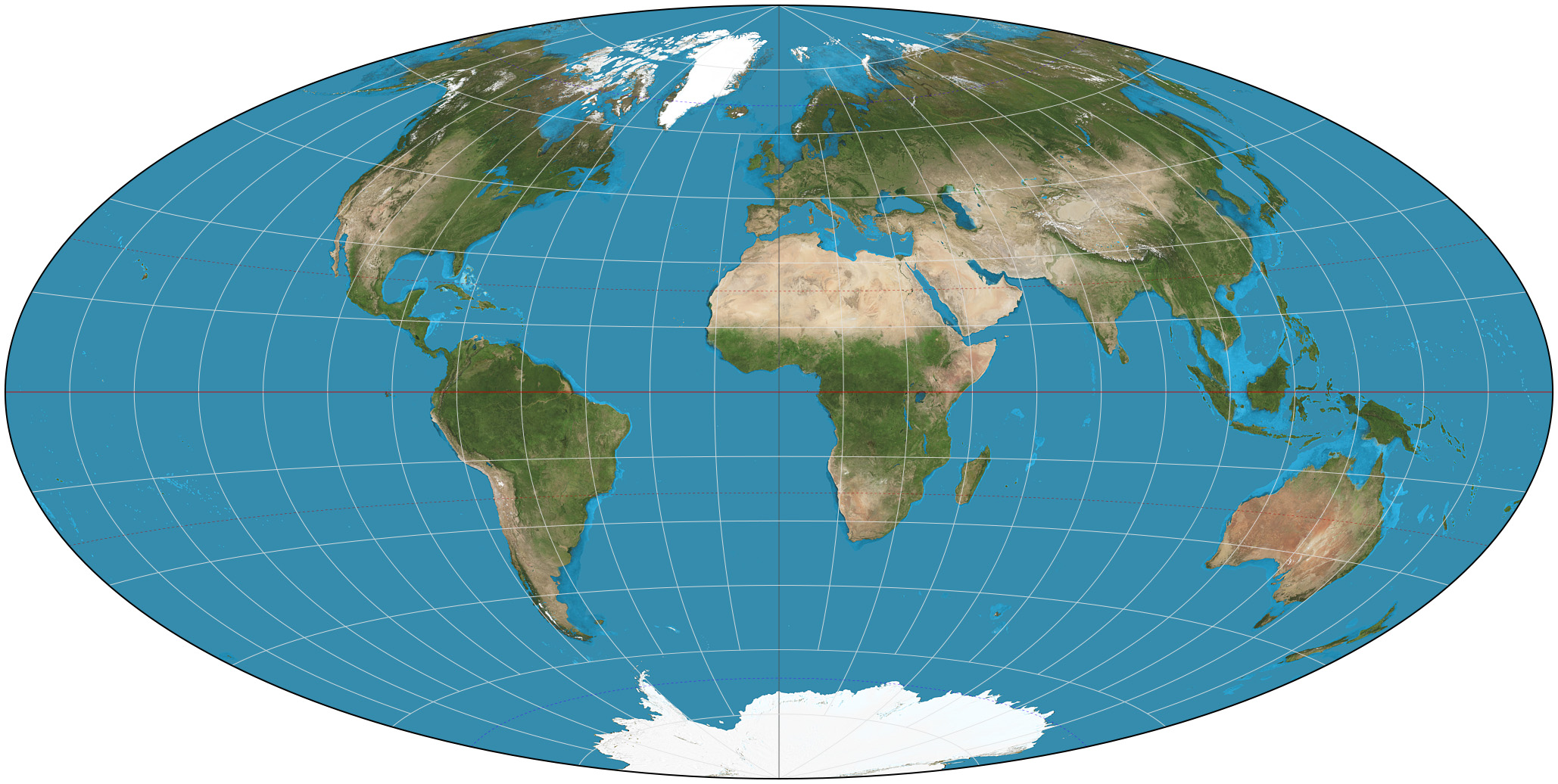} 
\\
\includegraphics[width=1.6in,align=c]{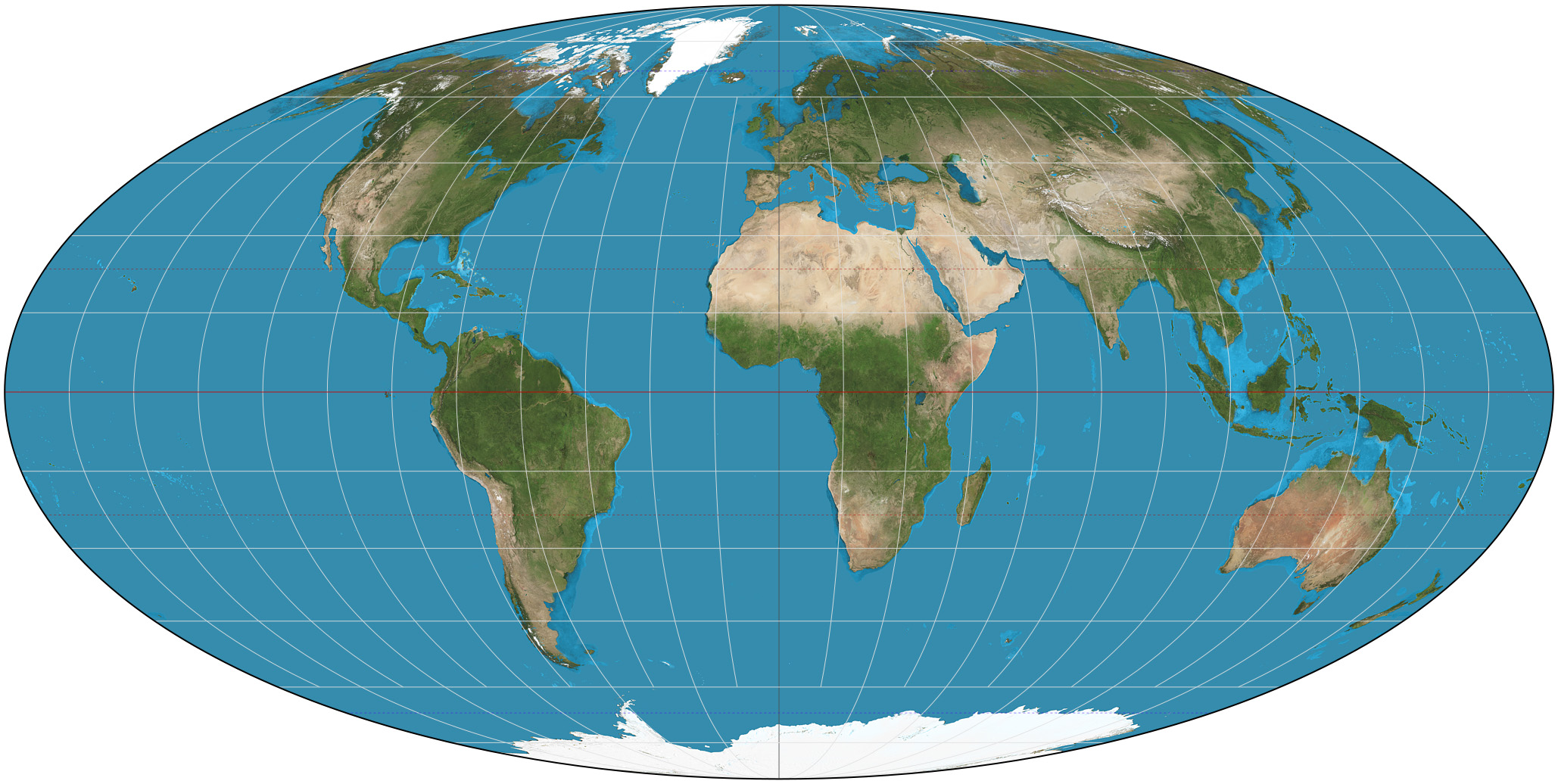}
\includegraphics[width=1.6in,align=c]{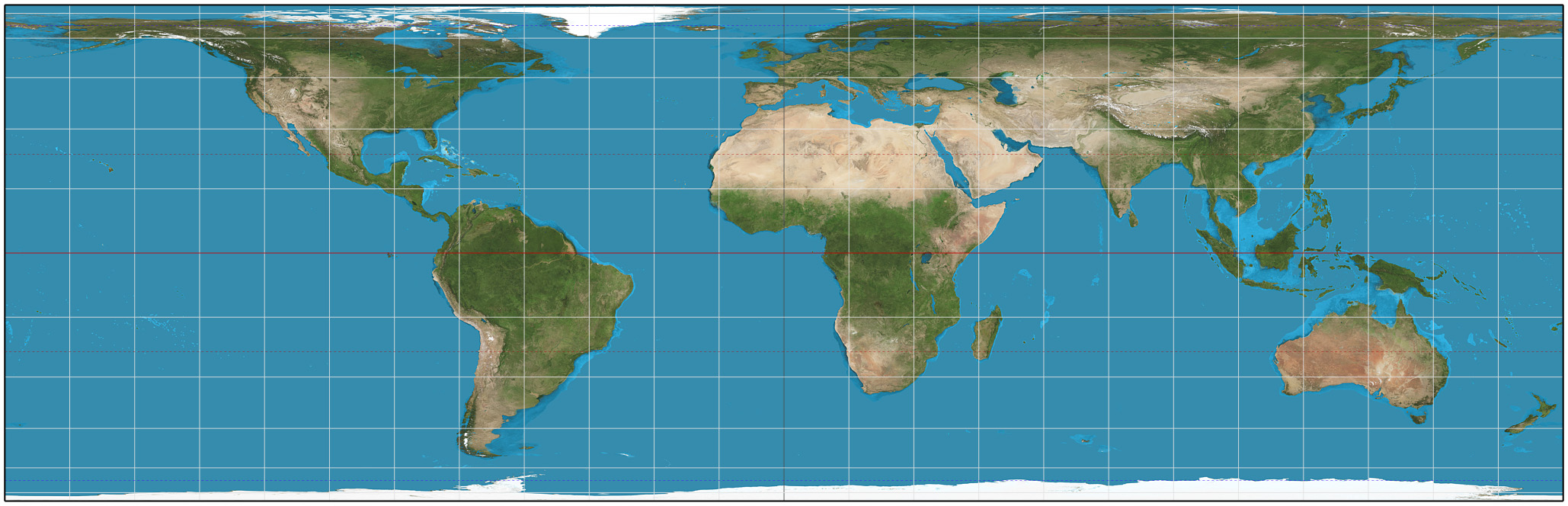}
\includegraphics[width=1.6in,align=c]{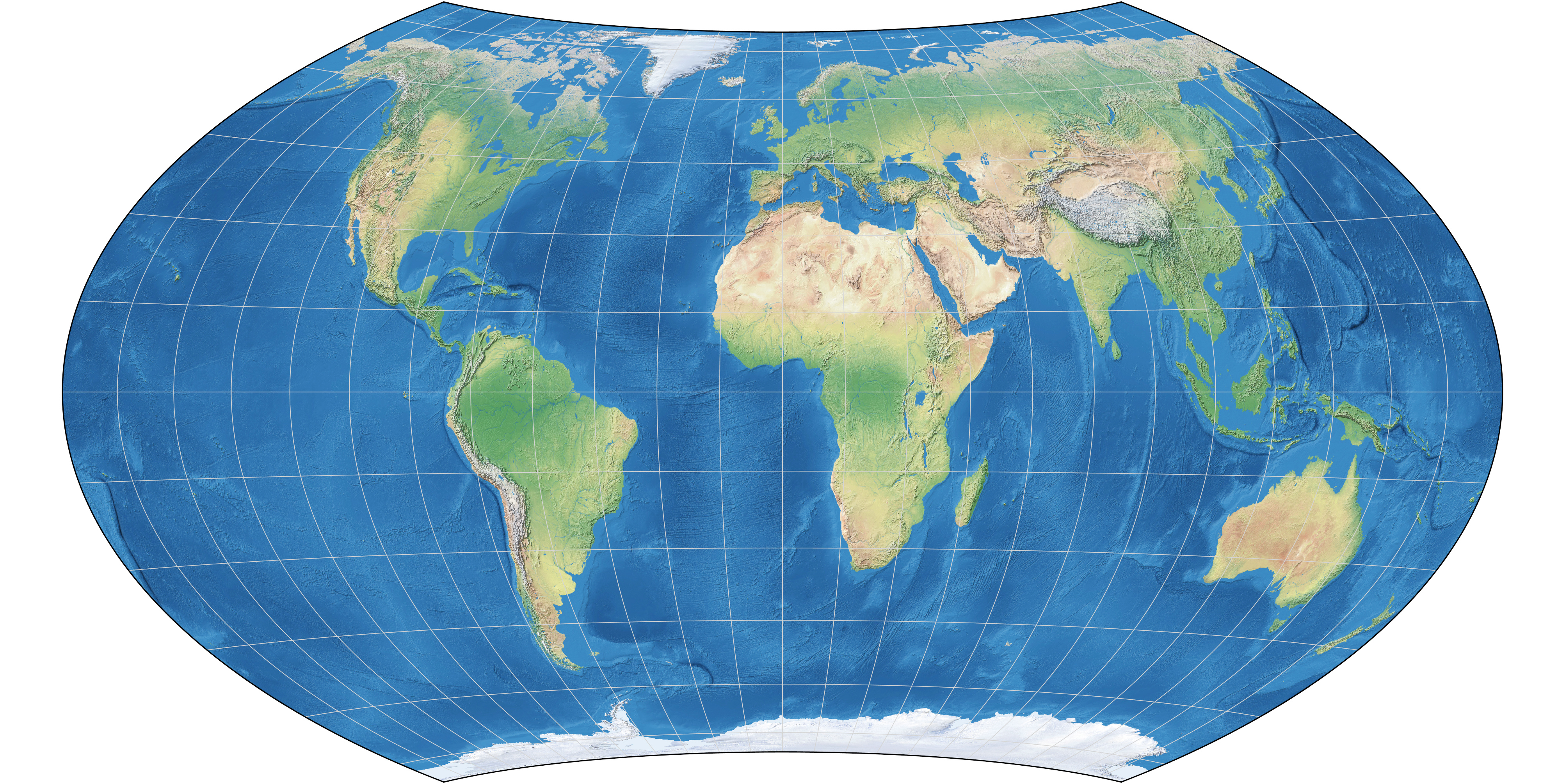}
\end{minipage}
\caption{Some map projections. Listed left-to-right, the top row shows the
        Equirectangular, Mercator and Guyou projections, the middle row shows
                the Winkel triple, Nicolosi and Aitoff projections, and the
                bottom row shows the Mollweide, Lambert and Wagner projections.
                \\ \phantom{Hi} \hfill {\em Image Credits: Daniel R. Strebe and Tom Patterson}
} \label{fig0}
\end{figure}

The fact that the Earth is spherical in shape makes it an interesting challenge
on how to make flat maps that accurately represent the not-flat Earth.
For maps of small areas, like cities or counties or states, it's easy to make a
flat map that gives an almost perfect representation of that surface.  But, how
can one make a map that shows the entire Earth on one flat map?  That is a
challenge.  Obviously no such map will be geometrically perfect.  Over hundreds
of years, several map projections have been proposed and analyzed for their quality.
And this work has created a research area called {\em cartography}.
Some of the most well-known projections are 
shown in Figure \ref{fig0}.
Until very recently, the Winkel tripel, shown as the left entry of the middle
row, was considered the ``best'' flat map of the Earth.

\section{A New Map Projection}

In \cite{GGV21}, a new flat map of the Earth (or any other spherical object) was
introduced that is superior to all other well-known flat maps as measured by
skewness, flexion, isotropy and area.  This
new map involves projecting one half of the sphere onto one side of a flat
circular disk using the azimuthal-equidistant projection and the other half onto
the other side using the same projection.  

We should note here that the azimuthal-equidistant projection is normally used to make a 
full-earth map.  If the north pole is positioned at the center of the map, then
the northern hemisphere occupies a central circular part of the map and the
southern hemisphere occupies an annulus that surrounds the northern hemisphere.
The south pole spans all around the outer edge of the map.   So, the distortion
at the south pole is infinitely large and therefore this map projection is not
considered one of the best.   But, the northern hemisphere part of the map is
very good and that is why the central half is used for the new two-sided flat map.

We should also note that if the two sides of the new two-sided map are placed
side-by-side to be viewed from one side, then it looks very similar to the
Nicolosi map shown in Figure \ref{fig0}.  But, they aren't identical.  The
differences can be found in the various cartography books that describe these
projections in detail.   See, for example \cite{cartography}.

To make the description of the new map
simple and easy to understand, it is helpful to assume that the 
equator is the separator of the two halves.  So, the northern hemisphere is
displayed flatly on one side and the southern hemisphere is displayed flatly on
the other side.  The longitude associated with each point along the edge is
the same when viewed from one side as it is on the other side.  Radial lines
from the center of the disk to the edge correspond to specific longitudes.  And,
the angular spacings of the longitude lines are equal.  Lastly, the
latitude varies linearly along the radial lines.  For example, the latitude
at the center is $\pm 90^\circ$, the latitude at the equatorial edge is
$0^\circ$ and the latitude halfway from the center to the edge is $\pm
45^\circ$.  Figure \ref{fig3} shows this flat disk map of the Earth.
If for each point on the sphere we let $\lambda$ denote the longitude and let 
$\theta$ denote $90^\circ$ minus the
latitude with both angles reexpressed in radians, then the GGV flattened map of the northern
hemisphere expressed in polar coordinates $(r, \phi)$ is this:
\[
    r = \theta,  \hspace*{0.3in}
    \phi = \lambda .
\]

This new map projection has generated a lot of interest in the cartography world and it was
featured in Time Magazine 
\cite{Time21}
as one of the 100 best inventions of 2021.

\section{A Stress Minimization Approach}

In this paper, we propose a small modification to the map projection described
above.
The modification is inspired by thinking about an interesting physical way in
which such a map might be made.  And the mathematics required to solve the problem
turns out to be a very interesting application of the calculus of variations and
an associated differential equation.   

Here's the physical approach.  Start with the
map printed on a spherical rubber ball.  Imagine that inside the rubber ball
there is a metal ring located at the equator and that this ring can be enlarged
as much as desired from its default size.  If the ring is just enlarged a little
bit then the ball will exhibit some oblateness.  As the ring gets enlarged more
and more, the northen and southern sides of the ball will start to get closer to
the equatorial plane.  At some point when the ring is sufficiently enlarged, the
ball will become perfectly flat.  The purpose of this paper is to do the physics
to determine exactly how far out from the center each latitudinal circle will
get stretched.


\begin{figure}
\includegraphics[width=5in]{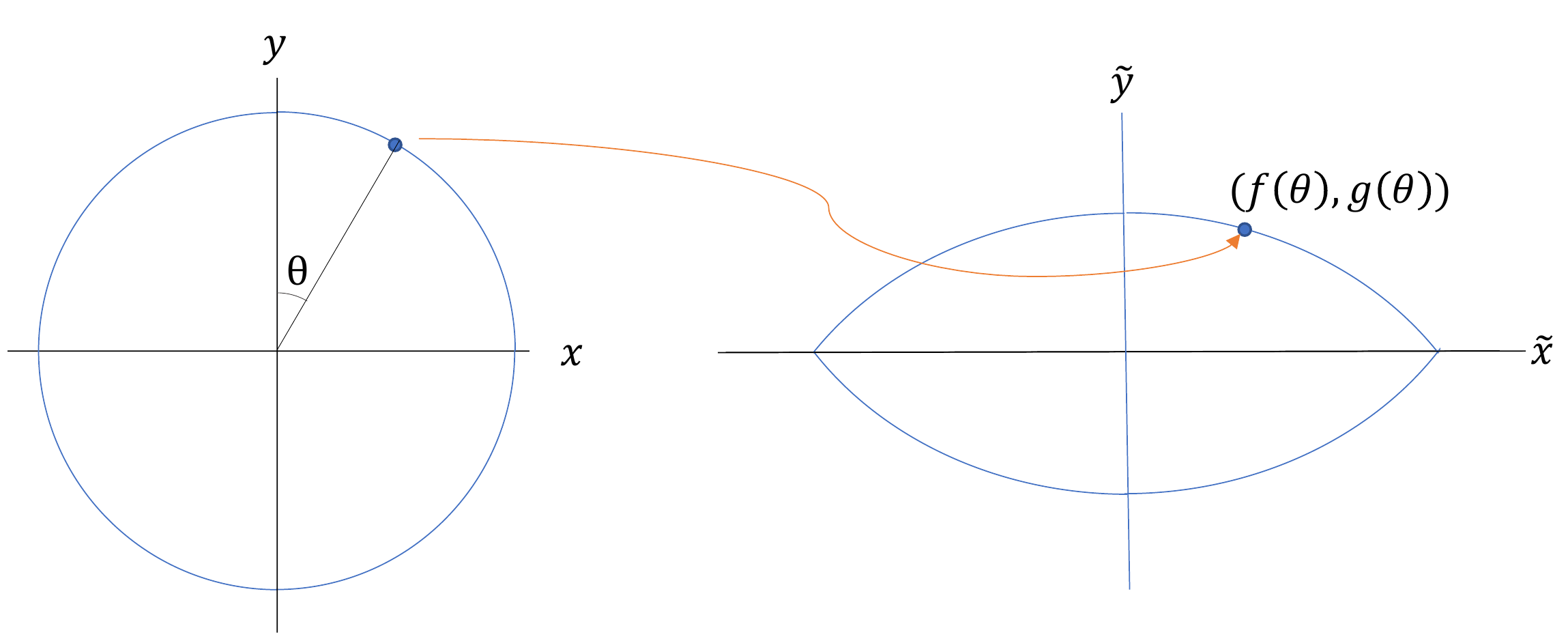}
\caption{A slice of the surface of the Earth rendered as a pale blue circle.   On the left
the Earth is shown before the equatorial stretch and on the right is
shown how it might appear after being stretched a bit.
In a fully flattened map, we would have $g(\theta) = 0$
for all $\theta$ in $[0, \pi/2]$.
} \label{fig2}
\end{figure}

Figure \ref{fig2} shows a single slice through the sphere.  The equatorial plane
is on the $x$-axis and the north/south poles are at the $y$-axis.  The blue
outline shows the surface of the Earth as it appears in this slice.  Our goal is
to determine where each point on the unstretched circle will appear on the
stretched out version.  So, let's consider the point at
\[
        x(\theta) = \sin(\theta),  \hspace*{0.3in} y(\theta) = \cos(\theta)
\]
and let's let $f$ and $g$ denote the corresponding location in the stretched circle:
\[
        \tilde{x}(\theta) = f(\theta), \hspace*{0.3in} \tilde{y}(\theta) = g(\theta) 
\]
According to physics, the shape of the stretched ball 
will be such that the integral
over the ball's surface of the magnitude squared of the stress tensor is minimized.
And, of course, it should be clear the $f$ is a smooth increasing function of $\theta$
and $g$ is a smooth decreasing function of $\theta$.   So, we will henceforth assume
that $f$ and $g$ are twice continuously differentiable ($f, g \in C^2[0,\pi/2]$), with $f(0) = 1$, $g(\pi/2) = 0$ and 
$f'(\theta) \ge 0$ and $g'(\theta) \le 0$ for $0 \le \theta \le \pi/2$.

At the point $(\tilde{x}(\theta), \tilde{y}(\theta))$ in the stretched circular slice, let
$\sigma(\theta)$ denote the stress in the direction tangent to the surface and let
$\rho(\theta)$ denote the stress in the direction perpendicular to the
$2$-dimensional plane of the slice.  We could also introduce notation for the
stress whose direction is in the plane of the slice but is perpendicular to the
tangential stress.  But we don't need this third component of the stress vector.
It is zero because the surface of the ball has no thickness and so there are no forces
perpendicular to the surface of the ball.

If we let $dx$ and $dy$ denote an infinitesimal section of the unstretched
circle corresponding to an infinitesimal angular segment of $d\theta$, then we have
\[
        \| (dx, dy) \| = \sqrt{dx^2 + dy^2} 
                       = \sqrt{\cos^2(\theta) + \sin^2(\theta)} \; d\theta 
                       = d \theta.
\]
Similarly, the corresponding length in the stretched circle is
\[
        \| (d\tilde{x}, d\tilde{y}) \| = \sqrt{ d\tilde{x}^2 + d\tilde{y}^2 } 
        = \sqrt{ f'(\theta)^2 + g'(\theta)^2 } \; d\theta .
\]
From these two displacements and assuming that the stress is linearly proportional to the displacement, the stress $\sigma(\theta)$ is
\[
        \sigma(\theta) = 
        \ds \frac{\| (d\tilde{x}, d\tilde{y}) \| }{ \| (dx, dy) \| } - 1
        = \sqrt{ f'(\theta)^2 + g'(\theta)^2 } \; - \; 1 .
\]
To determine the stress $\rho(\theta)$ perpendicular to the plane, note that the circumference of the circle on the sphere associated with the angle 
$\theta$ is $2 \pi \sin(\theta)$ 
and the corresponding circumference in the stretch sphere is 
$2 \pi f(\theta)$ and hence
\[
        \rho(\theta) = \frac{f(\theta)}{\sin(\theta)} \; - \; 1 .
\]
The total stress in the stretched ball is the integral over the entire sphere of
the norm squared of the stress vector which can easily be written as a single
integral
\[
    \int_0^{\pi/2} \left( \sigma(\theta)^2 + \rho(\theta)^2 \right) 2 \pi \sin(\theta) d \theta.
\]

Before we jump into the details of how to solve for the function $f$ that
minimizes this stress, let's make a few simple observations.  The first one is
that the function $\sin(\theta)$ appears in the denominator in the formula for $\rho(\theta)$.
This is a bit troubling because $\theta = 0$ is one of the values of $\theta$ in
the integral.  But, here's the thing... the numerator, $f(\theta)$, is also zero
at $\theta = 0$.   And, so we have $0/0$, which can be okay.  Or not.  It
depends on the slope of $f(\theta)$ as $\theta$ approaches zero.  If 
$f(\theta) \approx c \; \theta$ for some constant $c$, then the ratio does not
explode.   Of course, we could consider the case where $f(\theta) \approx c \;
\sqrt{\theta}$.  In that case, the function $\rho(\theta)$ goes to infinity as
$\theta$ tends to zero and hence the stress associated with such a choice of $f$
is infinite.   But, our aim is to minimize the stress.   And, at the minimum
this ratio is well behaved.
Henceforth, we will limit our analysis to functions $f$ (and perturbations
thereof) for which the stress integral is finite.

A second issue to discuss is the smoothness of the function $f$.  Since this is
a problem motivated by physics, it's probably safe to assume that $f$ is
infinitely differentiable.  But, in order to do the math, we only need to assume
that $f$ is twice continuously differentiable.   So, henceforth, that will be
our assumption and after we've found the minimal stress solution we can check to
see how smooth the function $f$ turns out to be.

\section{Solving the Stress Minimzation Problem}

\subsection{Calculus of Variations}

As shown above,
the total stress associated with the shape determined by the functions
$\sigma(\theta)$ and $\rho(\theta)$ is given by a simple integral.
In a fully flattened sphere, the function $g(\theta)$ is zero and so the stress
in that case, let's call it $S(f)$, is given by
\[
    S(f)
    = \int_0^{\pi/2} \left( \left( f'(\theta)-1 \right))^2 
                          + \left( \frac{f(\theta)}{\sin(\theta)} - 1 \right)^2 \right)
        2 \pi \sin(\theta) d \theta.
\]
Before we continue with the math, here's an important note... the function $g$
        that we introduced to represent the vertical component of the stretched ball's
        coordinates is now and henceforth zero and therefore we no longer need
        to think of $g$ as representing that component of the stretch and as we
        solve the minimum stress problem we will use $g$ to represent a
        completely different function that we introduce as we solve the
        differential equation.

Our goal now is to find the function $f(\theta)$ that minimizes the stress.
If $\theta$ were just a finite set of values between $0$ and $\pi/2$, then $f(\theta)$ could be 
thought of as a finite set of variables over which we wish to optimize the function $S(f)$.
In other words, if $\theta$'s domain was a finite set, then this problem could be solved using calculus (in 
rather high, but finite, dimension).  But, for our problem $\theta$ varies over
all real-valued numbers between 
$0$ and $\pi/2$.  In other words, this minimization problem involves an infinite
number (in fact a continuum) of variables.
Such problems can also be solved using the same approach that underlies calculus.  This generalized
methodology is called {\em calculus of variations}.  While we could just quote
the main formula that one needs to solve this kind problem, it's called the {\em
Euler-Lagrange equation}, we will not assume that this equation is
known and instead derive the solution in a logical manner that is
analogous to how calculus is used to solve optimization problems defined
over finite dimensional domains.
For readers who are comfortable with using the Euler-Lagrange equation, you can jump straight 
ahead to Equation \eqref{ode}.
%
%

Here's the derivation. Let $\Df(\theta)$, $0 \le \theta \le \pi/2$ denote a perturbation ``direction''.
Note:  in finite dimensional calculus, this would be a vector.
Given that the functions $f$ are assumed to be twice differentiable, the one that minimizes the stress function will be a function that's a
{\em critical point} of the stress function.  A critical point is one where an
infinitesimal perturbation in any direction is zero.
So, we need to find a function $f$ for which 
\[
        \lim_{\ve\rightarrow 0} \;
        \frac{S(f + \ve\Df) - S(f)}{\ve} 
        = 0 
\]
for all functions $\Df \in C^2[0,\pi/2]$ for which $\Df(0) = 0$.
There is this boundary condition on the function $\Df$ because the function $f$
must be zero at $\theta = 0$ and hence we are not allowed to perturb this value.
We will also restrict our attention to only those $f$ and $\Df$ for which $S(f)$ and
$S(f + \ve \Df)$ are finite.
Let's calculate this ratio:
\vspace*{0.1in}
\begin{eqnarray*}
        \setlength{\arraycolsep}{0.1em}
        \begin{array}{r c c c c l}
            \multicolumn{3}{l}{\ds \frac{S(f+\ve\Df)-S(f)}{\ve}} \\[0.2in]
            \hspace*{0.1in}
            = 2 \pi \stretchint{7ex}_0^{\pi/2} \Bigg(
                f'(\theta)^2 & + 2 f'(\theta) \ve\Df'(\theta) & +
                \ve^2\Df'(\theta)^2 & - 2 f'(\theta) & - 2 \ve\Df'(\theta) & + 1 \\[0.2in]
               -f'(\theta)^2 &                             &                  & + 2 f'(\theta) &                  & - 1 \\[0.2in]
               + \ds \frac{f(\theta)^2}{\sin^2(\theta)} &
               + 2 \ds \frac{f(\theta)\ve\Df(\theta)}{\sin^2(\theta)} &
               + \ds \frac{\ve^2\Df(\theta)^2}{\sin^2(\theta)} &
               - 2 \ds \frac{f(\theta)}{\sin(\theta)} &
               - 2 \ds \frac{\ve\Df(\theta)}{\sin(\theta)} &
               +1 \\[0.2in]
               - \ds \frac{f(\theta)^2}{\sin^2(\theta)} &
               &
               &
               + 2 \ds \frac{f(\theta)}{\sin(\theta)} &
               &
               -1 \Bigg) \sin(\theta) \ds \frac{1}{\ve} d \theta . \\[0.2in]
            \hspace*{0.1in}
            =  2 \pi \stretchint{7ex}_0^{\pi/2} \Bigg( \hspace*{0.4in}
              &  2 f'(\theta) \Df'(\theta) & +
                \ve\Df'(\theta)^2 & & - 2 \Df'(\theta) \\[0.2in]
              &   \ds 
               + 2 \ds \frac{f(\theta)\Df(\theta)}{\sin^2(\theta)} &
               + \ds \frac{\ve\Df(\theta)^2}{\sin^2(\theta)} & &
               - 2 \ds \frac{\Df(\theta)}{\sin(\theta)} &
               \hspace*{0.2in} \Bigg) \sin(\theta) \ds  d \theta . 
        \end{array}
        \\[0.2in] ~
\end{eqnarray*}
Now taking the limit as $\ve$ tends to zero, we get the following
expression for the perturbational change in the stress:
\vspace*{0.1in}
\begin{equation*}
        \lim_{\ve\rightarrow 0} \;
        \frac{S(f + \ve\Df) - S(f)}{\ve} 
    = 4 \pi \stretchint{7ex}_0^{\pi/2} \Bigg(
        f'(\theta) \Df'(\theta) - \Df'(\theta) 
       + \ds \frac{f(\theta)\Df(\theta)}{\sin^2(\theta)} 
       - \ds \frac{\Df(\theta)}{\sin(\theta)} 
       \Bigg) \sin(\theta) d \theta . 
\end{equation*}
\vspace*{0.2in}
Setting this perturbation to zero and dividing by $4 \pi$ we get
\begin{equation} \label{eq1}
  0 = \stretchint{7ex}_0^{\pi/2} 
        \left( f'(\theta) - 1 \right) \sin(\theta) \Df'(\theta) d\theta
       + 
       \stretchint{7ex}_0^{\pi/2} 
       \left( \ds \frac{f(\theta)}{\sin(\theta)} - 1 \right) \Df(\theta) d\theta 
       .
\end{equation}
Next we do integration by parts on the first integral to convert the
$\Df'(\theta)$ in the integrand to $\Df(\theta)$:
\begin{eqnarray*}
  \stretchint{5ex}_0^{\pi/2} 
  \left( f'(\theta) - 1 \right) \sin(\theta) \Df'(\theta) d\theta
  & = &
  \left( f'(\pi/2) - 1 \right) \Df(\pi/2)  \\
  &&
  \hspace*{-0.4in}
  - 
  \stretchint{5ex}_0^{\pi/2} 
  \left( f''(\theta) \sin(\theta) + (f'(\theta) - 1) \cos(\theta) \right) \Df(\theta) d\theta .
\end{eqnarray*}
Substituting this into Eq. \eqref{eq1}, we get
\begin{equation*} 
  0 = 
  \left( f'(\pi/2) - 1 \right) \Df(\pi/2)  
  - 
  \stretchint{7ex}_0^{\pi/2} 
  \left( 
          f''(\theta) \sin(\theta) + (f'(\theta) - 1) \cos(\theta) 
          - \ds \frac{f(\theta)}{\sin(\theta)} + 1 
  \right) \Df(\theta) d\theta .
\end{equation*}
For this expression to be zero for all 
$\Df \in C^2[0,\pi/2]$ for which $\Df(0) = 0$, 
it
must be true that the factors multiplying these perturbations are zero.
Hence, to find the 
        functions that are critical ``points'' of the stress function,
we need to solve this differential equation for the function $f$:
\vspace*{0.2in}
\begin{eqnarray} \label{ode}
    \sin^2(\theta) f''(\theta) + \sin(\theta)\cos(\theta) f'(\theta) - f(\theta)
    & = & \sin(\theta)\cos(\theta) - \sin(\theta) \\ \nonumber
    f(0) & = & 0 \\ \nonumber
    f'(\pi/2) & = & 1 .
\end{eqnarray}
\vspace*{0.2in}

\subsection{Solving the Differential Equation}

\hspace*{0.0in}
The differential equation \eqref{ode} is a second-order linear differential
        equation and hopefully not too difficult to solve.
But, the ``coefficients'' multiplying the terms involve sines and cosines and
therefore this differential equation might not be easy to solve.  Let's give it a try.

Before we actually solve the equation, let's note that we are expecting $f$ to
be close to, and maybe even equal to, a linear function with slope 1, 
   i.e., $f(\theta) = \theta$ for $0 \le \theta \le \pi/2$.  
With this function, $f'(\theta) = 1$ and $f''(\theta) = 0$ and so the left-hand
side of equation \eqref{ode} becomes $\sin(\theta) \cos(\theta) - \theta$ which
is interestingly similar to the right-hand side but not equal to it.

Using generic notation, a second-order linear differential equation can be
written like this:
\[
    \alpha_2(\theta) f''(\theta) + \alpha_1(\theta) f'(\theta) + \alpha_0(\theta)
            f(\theta) = \beta(\theta),
\]
where the functions $\alpha_i$, $i=0,1,2$, and $\beta$ are known
functions.
If $\alpha_2(\theta) = 0$ for all $\theta$ in the domain of interest, then the
problem would be a first-order differential equation and it would be fairly easy
to solve.  Unfortunately, there's no simple trick to zero out $\alpha_2$.
But, it is possible to zero out the $\alpha_0$ coefficient and that makes the
solution process easier.  We can get rid of $\alpha_0$ by rewriting the
differential equation in terms of a different function $g$ related to $f$ like
this
\[
    f(\theta) = \gamma(\theta) g(\theta)
\]
where $\gamma$ is an appropriately chosen function.  Let's write the
differential equation using the function $g$ and see what that tells us we need
to choose for the function $\gamma$.  Differentiating once we get
\[
    f' = \gamma g'  +   \gamma' g
\]
and differentiating a second time we get
\[
   f'' = \gamma g'' + 2 \gamma' g' + \gamma'' g .
\]
(Note that henceforth we will often not explicitly show the argument for
 functions of $\theta$.)
Substituting these expressions into \eqref{ode}, see that
\begin{eqnarray*}
    \sin^2(\theta) f'' + \sin(\theta)\cos(\theta) f' - f
    & = &
    \sin^2(\theta)             \left( \gamma g'' + 2 \gamma' g' + \gamma'' g \right) \\ && \hspace*{0.3in}
    + \sin(\theta)\cos(\theta) \left(                \gamma  g' + \gamma'  g \right) \\ && \hspace*{0.3in}
    -                                                          \; \gamma \;g \\
    & = &
               \sin^2(\theta)                                      \gamma                   g'' \\ && \hspace*{0.3in}
    + \left( 2 \sin^2(\theta) \gamma'  + \sin(\theta) \cos(\theta) \gamma         ) \right) g'  \\ && \hspace*{0.3in}
    + \left(   \sin^2(\theta) \gamma'' + \sin(\theta) \cos(\theta) \gamma' - \gamma \right) g .
\end{eqnarray*}
To get the coefficient multiplying $g$ to vanish, we need to find a function $\gamma$
that satisfies this differential equation:
\[
    \sin^2(\theta) \gamma'' + \sin(\theta)\cos(\theta) \gamma' - \gamma = 0 .
\]
In other words, we need to find a solution to the homogeneous variant of the
original differential equation.  It's not trivial, but it is possible.
Given that this equation only has sines and cosines in its coefficients, this
suggests that the function $\gamma$ is probably also a simple trigonometric function.
So, let's see if we can find a function $\beta$ for which
\[
        \gamma(\theta) = \beta(\sin(\theta))
\]
is a solution to the differential equation.
Taking first and second derivatives and plugging them into the differential equation for $\gamma$, we
get this version of the differential equation written using the function $\beta$:
\[
  \sin^2(\theta) \cos^2(\theta) \beta''(\sin(\theta)) 
  + \sin(\theta) \left( \cos^2(\theta) - \sin^2(\theta) \right) \beta'(\sin(\theta))
  - \beta(\sin(\theta)) = 0.
\]
And, if we let
\[
    x = \sin(\theta)
\]
we get this differential equation for $\beta$:
\begin{equation} \label{beta_ode}
    x^2(1-x^2) \beta''(x) + x(1-2 x^2) \beta'(x) - \beta(x) = 0.
\end{equation}
To find a solution to this differential equation, the normal next step would be to 
write the function $\beta$ as a power series 
\[
    \beta(x) = \sum_{j=0}^{\infty} a_j x^j
\]
and investigate the conditions that the coefficients must satisfy.
If we do this, we'll discover that all the even coefficients are zero, that
$a_1$ is anything, let's let it be $1$, and that $a_j = \frac{j-2}{j+1} a_{j-2}$
for all odd values of $j \ge 3$.
That's a rather complicated power series.  In fact, one can check that it's the
power series of this function:
\[
    \beta(x) = \frac{2x}{1+\sqrt{1-x^2}} .
\]
The function $\gamma$ associated with this function $\beta$ is
\[
    \gamma(\theta) = \frac{2\sin(\theta)}{1+\cos(\theta)} .
\]
We could move forward with this formula for $\gamma$.  But, it's a bit
complicated and one has to wonder why we didn't find two solutions and if we did
find a second solution would it be a simpler choice.  Let's answer these
questions.  It turns out that we didn't find a second independent solution
because it will have a singularity at $x=0$.  So, it doesn't have a power series
representation.  
Here's a nice way to find solutions that can be singular.  Instead of doing a
power series using a sum from zero to infinity, let's consider a sum from
$-\infty$ to $\infty$:
\[
    \beta(x) = \sum_{j=-\infty}^{\infty} a_j x^j.
\]
With this much broader sum we get these derivatives:
\[
    \beta'(x) = \sum_{j=-\infty}^{\infty} a_j j x^{j-1} 
    \hspace*{0.3in} \text{and} \hspace*{0.3in}
    \beta''(x) = \sum_{j=-\infty}^{\infty} a_j j(j-1) x^{j-2}.
\]
Plugging these into the differential equation \eqref{beta_ode} for $\beta$ and simplifying, we
get
\begin{eqnarray*}
    && x^2(1-x^2) \beta''(x) + x(1-2 x^2) \beta'(x) - \beta(x) \\ 
    && \hspace*{1in} = \sum_{j=-\infty}^{\infty} (j-1) \left( (j+1) a_j - (j-2) a_{j-2} \right) x^j \\ 
    && \hspace*{1in} = 0.
\end{eqnarray*}
The only way for this sum to be zero for all $x$ is for all of the coefficients
to be zero and so we get:
\[
    (j+1) a_j = (j-2) a_{j-2} \hspace*{0.2in} \text{for $j \ne 1$}.
\]
The solution we found before had $a_0 = 0$ and $a_1 = 1$.   If we put both of
these coefficients to zero, then it's easy to see that $a_j = 0$ for all 
$j \ge 0$.  Because there is no condition relating $a_1$ to $a_{-1}$, it follows
that we can set $a_{-1}$ to anything we like.   Let's set it to one.  It is now
interesting to note that $a_{-3}$ does not depend on $a_{-1}$ because the
equation relating them is $(-1+1)a_{-1} = (-1-3)a_{-3}$ which tells us that
$a_{-3} = 0$ no matter what value we assign to $a_{-1}$.  
And, recursing by twos to smaller and smaller values of $j$, we
see that $a_j = 0$ for all odd values of $j \le -3$.
Similarly, we can start with $a_0 = 0$ and recurse to smaller and smaller even
values of $j$ to discover that $a_j = 0$ for all even values of $j \le 0$.  So,
we have just discovered that 
\[
    \beta(x) = \frac{1}{x}
\]
is a solution to $\beta$'s differential equation.  And, hence, 
\[
    \gamma(\theta) = 1/\sin(\theta).
\]
This is much simpler than the other solution and so let's move forward with this
choice of $\gamma$.
With this choice we get this differential equation for $g$:
\[
    \sin(\theta) g'' 
    + \left( -2 \sin^2(\theta) \frac{\cos(\theta)}{\sin^2(\theta)} +
                    \sin(\theta)\cos(\theta)\frac{1}{\sin(\theta)} \right) g'
            = \sin(\theta) \cos(\theta) - \sin(\theta).
\]
And, it's easy to check that the coefficient multiplying $g'$ simplifies nicely to just
$-\cos(\theta)$.
%
Hence, the function $f$ satisfies its differential equation 
if and only if $g$ solves this differential equation:
\begin{eqnarray*}
    \sin(\theta) g''(\theta) - \cos(\theta) g'(\theta) 
    & = & \sin(\theta)\cos(\theta) - \sin(\theta) \\
    g(0) & = & 0 \\
    g'(\pi/2) & = & 1 .
\end{eqnarray*}
Now, to solve this differential equation, let
\[
    h(\theta) = g'(\theta).
\]
The function $h$ is the solution to this differential equation:
\begin{eqnarray*}
    \sin(\theta) h'(\theta) - \cos(\theta) h(\theta) 
    & = & \sin(\theta)\cos(\theta) - \sin(\theta) \\
    h(\pi/2) & = & 1.
\end{eqnarray*}
If we divide both sides of this differential equation by $\sin^2(\theta)$, then on the left we have
the derivative of $h(\theta)/\sin(\theta)$ and therefore this ratio can be determined by simply 
integrating the explicit function of $\theta$ on the right.
Doing this integration, we get that the
solution to this differential equation is:
\[
    h(\theta) = \sin(\theta) + \log(\cos(\theta) + 1) \sin(\theta) .
\]
Now that we have $h$, we integrate it to get the formula for $g$:
\begin{eqnarray*}
    g(\theta) & = & \int_0^{\theta} g'(\theta) d\theta + g(0) \\
              & = & \int_0^{\theta} h(\theta) d\theta \\
              & = & \int_0^{\theta} \left(
                              \sin(\theta) + \log(\cos(\theta) + 1) \sin(\theta) 
                              \right) d\theta \\
              & = & - \cos(\theta) + 1
                    + \int_0^{\theta} 
                              \log(\cos(\theta) + 1) \sin(\theta) 
                              d\theta .
\end{eqnarray*}
To do this last integral, let 
\[
    u = \cos(\theta) + 1.
\]
Then, $du = -\sin(\theta) d\theta$ and therefore
\begin{eqnarray*}
    \int_0^{\theta} 
              \log(\cos(\theta) + 1) \sin(\theta) 
              d\theta 
    & = & - \int_2^{\cos(\theta)+1} \log(u) du \\
    & = & - \Big[ u \log(u) - u \Big]_2^{\cos(\theta)+1} \\[0.1in]
    & = & - (\cos(\theta) + 1) \log( \cos(\theta) + 1 ) + \cos(\theta) - 1 + 2 \log(2).
\end{eqnarray*}
\vspace*{0.2in}
Hence,
\[
    g(\theta) = 2 \log(2) - \left( \cos(\theta) + 1 \right) \log( \cos(\theta) + 1 )
\]
and from this we get
\[
    f(\theta) = \frac{2 \log(2) - \left( \cos(\theta) + 1 \right) \log( \cos(\theta) + 1 )}{\sin(\theta)}.
\]
We now have an explicit solution for the one, and only, critical point for the
        stress function.  
        Given that the stress function is an integral of a positively weighted
        sum of squares of linear functionals 
        of the unknown function $f$, it seems pretty obvious that this critical
        point is the global minimum.   But, in case there are doubts, let's show that
        the ``second derivative'' is
        nonnegative no matter what ``direction'' perturbation we choose.  In
        other words, we will show that for every choice of the functions $f$ and $\Df$ we have
        \vspace*{0.1in}
        \[
            \lim_{\ve\rightarrow 0} \;\;
            \ds \frac{S(f+\ve\Df)-2S(f)+S(f-\ve\Df)}{\ve^2}
            \;\; \ge \;\; 0.
        \]
        To this end, let's do some algebra:
        \vspace*{0.1in}
        \begin{eqnarray*}
                \setlength{\arraycolsep}{0.1em}
                \begin{array}{r c c c c l}
                    \multicolumn{3}{l}{
                            \ds \frac{S(f+\ve\Df)-2S(f)+S(f-\ve\Df)}{\ve^2}
                    } \\[0.2in]
                    \multicolumn{3}{l}{
                    \hspace*{0.3in}
                    =  \ds \frac{S(f+\ve\Df)-S(f)}{\ve^2}
                       \; + \;
                       \ds \frac{S(f-\ve\Df)-S(f)}{\ve^2}
                    } \\[0.2in]
                    \hspace*{0.3in}
                    = \;\;\; \ds \frac{2 \pi}{\ve^2} \stretchint{7ex}_0^{\pi/2} \Bigg( 
                      &  \; 2 f'(\theta) \ve\Df'(\theta) & + \;
                        \ve^2\Df'(\theta)^2 & & - \; 2 \ve\Df'(\theta) \\[0.2in]
                      &   \ds 
                       + \; 2 \ds \frac{f(\theta)\ve\Df(\theta)}{\sin^2(\theta)} &
                       + \; \ds \frac{\ve^2\Df(\theta)^2}{\sin^2(\theta)} & &
                       - \; 2 \ds \frac{\ve\Df(\theta)}{\sin(\theta)} &
                       \Bigg) \sin(\theta) \ds  d \theta \\[0.2in]
                    \hspace*{0.3in}
                    \;\;\; + \;  \ds \frac{2 \pi}{\ve^2} \stretchint{7ex}_0^{\pi/2} \Bigg( 
                      &- \; 2 f'(\theta) \ve\Df'(\theta) & + \;
                        \ve^2\Df'(\theta)^2 & & + \; 2 \ve\Df'(\theta) \\[0.2in]
                      &   \ds 
                       - \; 2 \ds \frac{f(\theta)\ve\Df(\theta)}{\sin^2(\theta)} &
                       + \; \ds \frac{\ve^2\Df(\theta)^2}{\sin^2(\theta)} & &
                       + \; 2 \ds \frac{\ve\Df(\theta)}{\sin(\theta)} &
                       \Bigg) \sin(\theta) \ds  d \theta \\[0.2in]
                    \multicolumn{5}{l}{
                    \hspace*{0.3in}
                    =  4 \pi \stretchint{7ex}_0^{\pi/2} \Bigg( 
                        \Df'(\theta)^2 + \;  
                        \ds \frac{\Df(\theta)^2}{\sin^2(\theta)} 
                       \Bigg) \sin(\theta) \ds  d \theta 
                    } \\[0.3in]
                    \multicolumn{1}{l}{
                    \hspace*{0.3in}
                    \ge 0 .
                    }
                \end{array}
                \\[0.2in] ~
        \end{eqnarray*}
        Because the function $S$ is an integral of squares of linear
        functionals, this second order differential turns out to be
        independent of the size $\ve$ of the perturbation.
        And, it is clear that the integrand is nonnegative and therefore so is
        the integral.  Hence, the critical point that we have found is indeed
        the global minimum.

\begin{figure}
\includegraphics[width=5in]{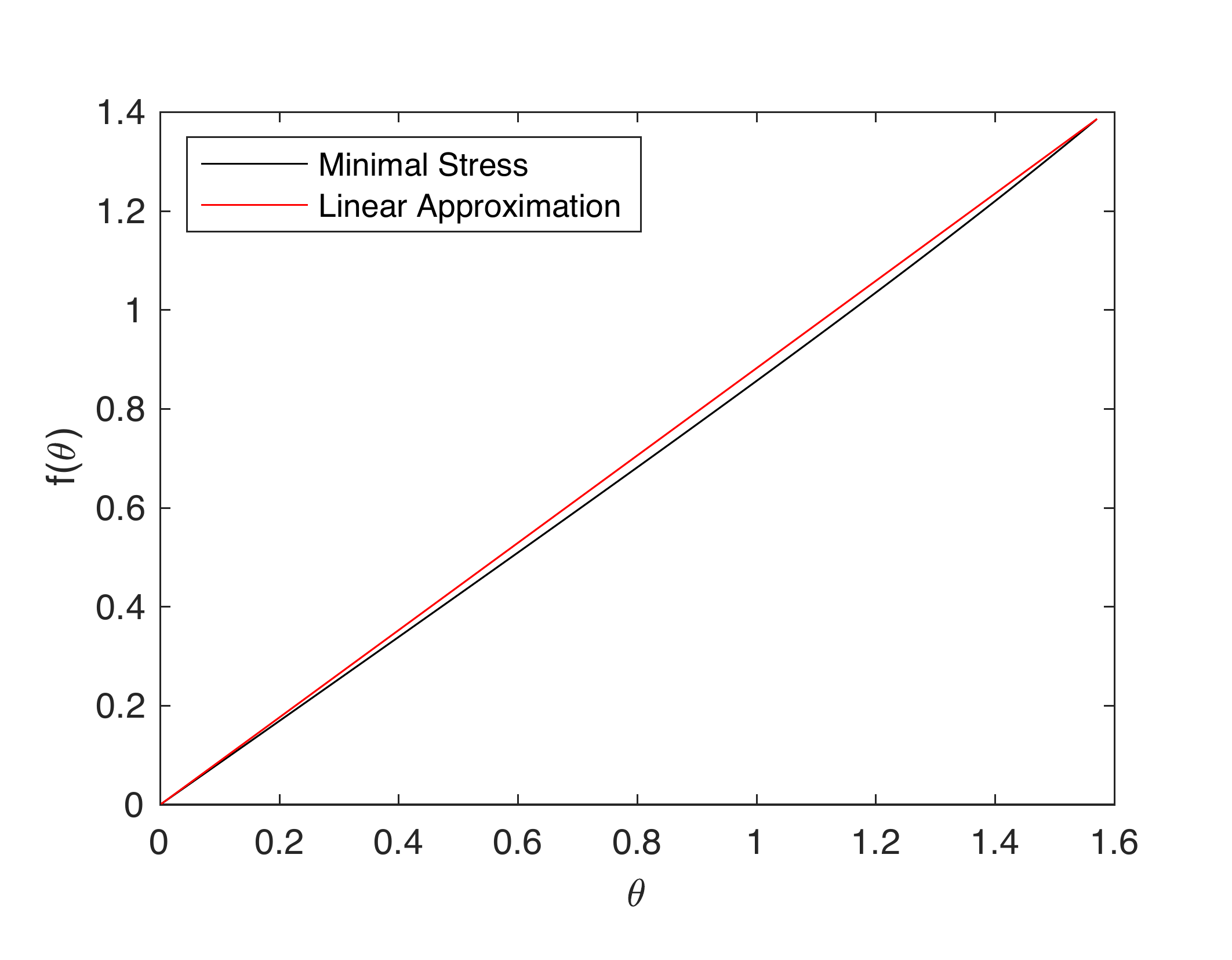}
\caption{The function $f(\theta)$ is shown in black.  It is almost a perfectly
        straight line with slope a little bit less than one.   But, it's not
        exactly a straight line.  The red line shown is a straight line
        that connects to the two end points of the function $f$.}
        \label{fig1}
\end{figure}

The function $f$ is shown in Figure \ref{fig1}.
The straight line, shown in red, is what one would use to do the
Gott-Goldberg-Vanderbei projection.  It is interesting how close these two curves are
to each other but yet are not the same.

It is also interesting to note that $f(\pi/2) = 2 \log(2) \approx 1.3868$ and hence the rubber ball
of radius $1$ got stretched to a disk having radius $2 \log(2)$.
Of course, we could have started with a smaller ball.   For example, had the original ball
had radius $1/(2 \log(2))$, then the flattened ball would have radius $1$.
So, the scale of the transformation doesn't matter.   What's interesting is how
close this function is to a linear function.

\noindent

\subsection*{Two Final Notes}  
For those readers who have access to {\it Mathematica}, the differential
equation for $f$ can be solved in just two lines of code..

\vspace*{0.2in}
\Small
\begin{verbatim}
    s = DSolve[ {Sin[x]^2*y''[x]+Sin[x]*Cos[x]*y'[x]-y[x]==Sin[x]*Cos[x]-Sin[x],
                           y[0]==0, y'[Pi/2]==1}, y[x], x] // FullSimplify
    f[x_]=y[x]/.s[[1]]
\end{verbatim}
\normalsize
\vspace*{0.1in}

\noindent
The output produced by {\it Mathematica} (with $x$ changed to $\theta$) is
\[
        f(\theta) =
        \log(2) \tan(\theta/2) 
        - 2 \cot(\theta/2)
        \log( \cos(\theta/2) ) .   
\]
It is easy to check that this solution is the same as the one derived above.

And, for those readers who have access to {\it Matlab}, the differential equation for $f$ can be also be solved
with just a few lines of code..

\vspace*{0.2in}
\Small
\begin{verbatim}
    syms f(x)
    f1 = diff(f,x);
    f2 = diff(f,x,2);
    ode = sin(x)^2 * f2 + sin(x)*cos(x) * f1 - f == sin(x)*cos(x) - sin(x);
    cond1 = f(0) == 0;
    cond2 = f1(pi/2) == 1;
    conds = [cond1 cond2];
    fSol(x)  = dsolve(ode,conds)
    fSim(x)  = simplify(fSol(x), 'steps', 14)
\end{verbatim}
\normalsize
\vspace*{0.1in}

\noindent
The output produced by {\it Matlab} (again with $x$ changed to $\theta$) is
%
\[
        f(\theta) = -\frac{\log(\cos(\theta)/4+1/4) + 2\log(e^{i \theta}+1)\cos(\theta) - \log(2)\cos(\theta) - \theta\cos(\theta)i}{\sin(\theta)}
\]
It is a bit of a challenge but it is also possible to show that this solution is the same as the one derived above.

\noindent

\subsection*{Acknowledgement}  
I would like to thank J. Richard Gott for many discussions which inspired me to give some serious thought to finding the best possible map projections.

\vfill
\pagebreak
\bibliographystyle{siam}
\bibliography{refs}

\end{document}